\documentclass{article}\begin{document}
\centerline{\bf Algebraic Equations Solved with Jacobi Elliptic Functions}
\vskip .4in

\centerline{Nikos Bagis}

\centerline{Stenimahou 5 Edessa}
\centerline{Pella 58200, Greece}
\centerline{nikosbagis@hotmail.gr}
\vskip .2in

\[
\]
\textbf{Keywords}: Elliptic functions; Quintic equation; Singular modulus; Modular equations;
\[
\]
\centerline{\bf Abstract}

\begin{quote}
In this article we solve a class of two parameter polynomial-quintic equation. The solution follows if we consider the Jacobian elliptic function $sn$ and relate it with the coefficients of the equation. The solution is the elliptic singular modulus $k$.   
\end{quote}

\section{Introduction and Definitions}
We will solve the following quintic equation  
\begin{equation}
e+h d Y+c Y^2+h b Y^3+aY^4+h Y^5=0 
\end{equation}
with respect to $Y$, where
\begin{equation}
e=\sqrt{1-x^2} \left(1-4 x^2+6 x^4-4 x^6+x^8\right)
\end{equation}
\begin{equation}
d=1-6 x^4+8 x^6-3 x^8\textrm{ , }
c=\sqrt{1-x^2} \left(-2+8 x^2-6 x^4\right)
\end{equation}
\begin{equation}
b=-2+6 x^4-4 x^6\textrm{ , }a=(1-4x^2) \sqrt{1-x^2} .
\end{equation}
Using the theory of Jacobian elliptic functions $sn, cn, dn$ and their addition property we construct the solution.\\ Our method is similar to that of solving equations using the trigonometric functions $\cos$ and $\sin$. For example, for the trisection of an angle $\theta$ holds the following formula
\begin{equation}
\cos(3\theta)=-3\cos(\theta)+4\cos^3(\theta) ,
\end{equation}
which rise from the addition formulas for the $\cos(x+y)$ and $\sin(x+y)$.\\
In Jacobian elliptic function theory we have the addition formulas\\
\begin{equation}
sn(u_1+u_2)=\frac{sn(u_1)cn(u_2)dn(u_2)+sn(u_2)cn(u_1)dn(u_1)}{1-k^2sn^2(u_1)sn^2(u_2)}
\end{equation}
\begin{equation}
cn(u_1+u_2)=\frac{cn(u_1)cn(u_2)-sn(u_1)sn(u_2)dn(u_1)dn(u_2)}{1-k^2sn^2(u_1)sn^2(u_2)}
\end{equation}
\begin{equation}
dn(u_1+u_2)=\frac{dn(u_1)dn(u_2)-k^2sn(u_1)sn(u_2)cn(u_1)cn(u_2)}{1-k^2sn^2(u_1)sn^2(u_2)}
\end{equation} 
where $k$ is the singular modulus. This defined as 
\begin{equation}
k^2=k_r^2=m=m(q)=\textrm{InverseEllipticNomeQ}[q]\textrm{, }q=e^{-\pi\sqrt{r}}\textrm{, }r>0
\end{equation}
and is solution of the equation\\
\begin{equation}
\frac{K(1-x)}{K(x)}=\sqrt{r},\textrm{ }r>0.
\end{equation} 
Note that if $r$ is positive rational then $k_r$ is algebraic number.\\
Also
\begin{equation}
u=\int^{am(u)}_{0}\frac{d\theta}{\sqrt{1-k^2\sin(\theta)^2}}
\end{equation}
\begin{equation}
sn(u)=\sin(am(u))\textrm{, }cn(u)=\cos(am(u))\textrm{, }dn(u)=\frac{d}{du}(am(u))
\end{equation}
\begin{equation}
sn(u,m)=\frac{2\pi}{Kk}\sum^{\infty}_{n=0}\frac{q^{n+1/2}\sin((2n+1)z)}{1-q^{2n+1}}
\end{equation}
\begin{equation}
cn(u)=\frac{2\pi}{Kk}\sum^{\infty}_{n=0}\frac{q^{n+1/2}\cos((2n+1)z)}{1+q^{2n-1}}
\end{equation}
\begin{equation}
dn(u)=\frac{\pi}{2K}+\frac{2\pi}{K}\sum^{\infty}_{n=1}\frac{q^{n}\cos(2nz)}{1+q^{2n}}
\end{equation}
where $u=2Kz/\pi$, $K=K(m)=K(k^2)$.
\begin{equation}
K(x)=\frac{\pi}{2}{}_{2}F_{1}\left[\frac{1}{2},\frac{1}{2};1;x\right],
\end{equation} 
where ${}_2F_{1}[a,b;c;z]$ is the Gauss hypergeometric function.

\section{Theorems}

We set $dn(u)=x$, $cn(u)=\alpha$, $sn(u)=\beta$, then hold (see [2]) $\alpha^2+\beta^2=1$ and $k^2\beta^2+x^2=1$. Hence in general 
\begin{equation}
\alpha=\frac{\sqrt{k^2+x^2-1}}{k}\textrm{ and }\beta=\frac{\sqrt{1-x^2}}{k}
\end{equation}
From the addition identities we have
\begin{equation}
sn(2u)=\frac{2 x \sqrt{1-x^2} \sqrt{k^2+x^2-1}}{k^2-\left(1-x^2\right)^2}
\end{equation}
\begin{equation}
cn(2u)=\frac{k^2+x^4-1}{k^2-\left(1-x^2\right)^2}
\end{equation}
\begin{equation}
dn(2u)=\frac{\left(1-x^2\right)^2-k^2 \left(1-2 x^2\right)}{k^2-\left(1-x^2\right)^2}
\end{equation}
and
\begin{equation}
sn(3u)=\frac{\sqrt{1-x^2}}{k}\frac{k^4 \left(1-4 x^2\right)+\left(1-x^2\right)^4+k^2 \left(-2+8 x^2-6 x^4\right)}{-\left(1-k^2\right)^2+6 \left(1-k^2\right) x^4-4 \left(2-k^2\right) x^6+3 x^8}
\end{equation}
\begin{equation}
cn(3u)=\frac{\sqrt{-1+k^2+x^2} \left(k^4+\left(-1+x^2\right)^4+k^2 \left(-2+4 x^2-6 x^4+4 x^6\right)\right)}{k \left(\left(-1+k^2\right)^2+6 \left(-1+k^2\right) x^4-4 \left(-2+k^2\right) x^6-3 x^8\right)}
\end{equation}
\begin{equation}
dn(3u)=\frac{x \left(6 k^2 \left(-1+x^2\right)^2+\left(-1+x^2\right)^3 \left(3+x^2\right)+k^4 \left(-3+4 x^2\right)\right)}{\left(-1+k^2\right)^2+6 \left(-1+k^2\right) x^4+8 x^6-7 x^8}
\end{equation}
Setting $h=sn(3u)$ we lead to equation (1). Hence the next theorem is valid:\\
\\
\textbf{Main Theorem.}\\
If $h=sn(3u)$ and $x=dn(u)$ are parameters related with the coefficients of  
$$
e+hdY+cY^2+hb Y^3+aY^4+hY^5=0 \eqno{(eq)}
$$
then a root of $(eq)$ is $Y=k=k_r$. Hence $(eq)$ is equivalent to
\begin{equation}
\frac{sn^{(-1)}(h,m_0)}{dn^{(-1)}(x,m_0)}=3\textrm{ and }Y=\sqrt{m_0}
\end{equation}
The $a,b,c,d,e$ are that of relations (2),(3),(4).\\
\\
\textbf{Remark.}\\
It holds from (10) that equation (24) is equivalent to
\begin{equation}
dn\left(\frac{1}{3}F\left[\arcsin(h),m_0\right],m_0\right)=x
\end{equation}
and
\begin{equation}
F\left[\arcsin\left(\sqrt{\frac{1-x^2}{m_0}}\right),m_0\right]=\frac{1}{3}F\left[\arcsin\left(h\right),m_0\right]
\end{equation}
where $F$ is the incomplete elliptic integral of the first kind.\\
\\
\textbf{Proposition 1.}\\
The value of $dn\left(\frac{1}{3}K,m\right)=dn\left(\frac{1}{3}K,k^2\right)$ is $L(k)$, where $L$ is known algebraic function.\\
\textbf{Proof.}\\
Set $h=sn(3u,k_r^2)=1$ and $Y=k_r$, then (1) is solvable with respect to $x$ in radicals. We evaluate $x=dn(u,k_r^2)=L(k)$. Using now (25) with $h=1$ we get $u=K$ and the result follows.\\ 
\\
\textbf{Note.} The value of $L$ is quite complicated in radicals but one can found it always. The fact that $dn(\frac{1}{3}K,m)=dn(\frac{1}{3}K,k^2)$ is $L(k)$ increases the complexity of the problem, but always keeping it solvable, since is function of the singular modulus.\\
In general for all positive real $r$ we have
$$
dn^2\left(\frac{1}{3}K,k^2\right)=1-k^2-\sqrt{\frac{(kk')^{4/3}}{2^{2/3}}+k^4-k^2}+
$$
\begin{equation}
+\frac{1}{2} \sqrt{-2 \sqrt[3]{2} (kk')^{4/3}+\frac{4 \left(2 k^6-3 k^4+k^2\right)}{\sqrt{\frac{(k k')^{4/3}}{2^{2/3}}+k^4-k^2}}+8 k'^4-8 k'^2}
\end{equation}
\\   
\textbf{Proposition 2.}\\
If $dn\left(\frac{1}{3}K,k_r^2\right)$ is rational then we can find always $k_r$.\\
\\
\textbf{Proposition 3.}\\
Let $\lambda=\frac{dn(u)}{dn(3u)}$ is known. In the case which $k_r$ is also known, then we can find $dn(u)$.\\
\textbf{Proof.}\\
Set $x=\lambda dn(3u)$ in (21) then one can see (for example with Mathematica program) that the reduced equation is solvable in radicals with respect to $dn(3u)$ and holds $dn(3u)=\Phi(\lambda,k_{r})$, where $\Phi$ known. Hence we find $dn(3u)$ and also $dn(u)=x$. The values of $a$ and $b$ are also follow from (17).\\ 
\\ 
\textbf{Proposition 4.}\\
If we know $dn(u)$ and $dn(3u)$, then we can find $k$ in radicals.\\
\textbf{Proof.}\\
Observe that knowing $dn(3u)$ and $dn(u)$, equation (24) is again solvable with respect to $k$.\\
\\
\textbf{Proposition 5.}\\
If we know $dn(u)$ and $dn(4u)$, then we can find $k$ in radicals.\\
\\
\textbf{Examples.}\\
1) The equation 
\begin{equation}
243\sqrt{3}-756hY-720\sqrt{3}Y^2+1728hY^3+ 
512\sqrt{3}Y^4-1024hY^5=0 
\end{equation}
with $h=-\frac{33 \sqrt{\frac{3}{2}}}{37}$ is solvable with elliptic functions. Comparing the coefficients with (2),(3),(4) we find
$$
dn(u)=x=1/2\textrm{ together with }h=sn(3u)=-\frac{33\sqrt{\frac{3}{2}}}{37}
$$ 
hence $Y=k=\frac{1}{\sqrt{2}}$ and $q=e^{-\pi}$.\\From this we have $dn\left(u,\frac{1}{2}\right)=\frac{1}{2}$ and $
sn\left(3u,\frac{1}{2}\right)=-\frac{33\sqrt{\frac{3}{2}}}{37}
$\\
\\
2) The equation 
$$
\sqrt{2} \left(-265+153 \sqrt{3}\right)+\left(500-288 \sqrt{3}\right) Y-32 \sqrt{2} \left(-17+9 \sqrt{3}\right) Y^2+
$$
$$
+64 \left(-16+9 \sqrt{3}\right) Y^3-512 \sqrt{2} Y^4+1024 Y^5=0
$$
have $h=1$ and $x=\frac{\sqrt{2+\sqrt{3}}}{2}$, hence $Y=k_{5/3}=\frac{1}{4} \sqrt{8-\sqrt{\frac{3}{2} \left(27-7 \sqrt{5}\right)}}$\\
3) For $h=1$ and $x=\sqrt{-2449+343 \sqrt{51}+20 \sqrt{6 \left(4999-700 \sqrt{51}\right)}}$, we get $Y=k_{34/3}$.\\
One can find $k_{34/3}$ with the program Mathematica. Taking the numerical value $k_{34/3}$ with 1000 digits accuracy and then using the command 'Recognize' we lead to a symmetric octic equation of $Y$. Mathematica do not solve this equation but with the change of variable $Y=iW$, ($i=\sqrt{-1}$) and then $W=v+v^{-1}$, the equation becomes solvable after calculating the powers $v^{\nu}+v^{-\nu}$ as polynomials of $v+v^{-1}$.\\
\\
\textbf{Proposition 6.}\\
It is $sd(u)=\frac{sn(u)}{dn(u)}$ and suppose that  
$sd(3u)=1-9x^2+24x^4-16x^6$, where $x=dn(u)$, $x$ known. Then we can find always $k=k_r$.\\
\\
\textbf{Theorem 1.}\\ 
The equation
$$
h k[k^4+k^2 \left(-4 X^3+6 X^2-2\right)+(1-X)^3 (3 X+1)]+[k^4 (1-4 X)+
$$
\begin{equation}
+k^2 \left(-6 X^2+8 X-2\right)+(1-X)^4]\sqrt{1-X}=0
\end{equation}   
have solution with respect to $X$
\begin{equation}
X=dn^2\left(\frac{1}{3}F\left[\arcsin(h),k^2\right],k^2\right)
\end{equation}   
\\
\textbf{Example.}\\
With $h=1/2$ equation (29) becomes
$$
2 X^{15}+18 X^{13}+42 X^{11}-38 X^9+7 X^8-186 X^7-8 X^6+X^5 \left(12 k^2+54\right)+
$$
\begin{equation}
+X^4 \left(6-6 k^2\right)+X^3 \left(8 k^4-24 k^2+270\right)+X\left(-6 k^4+12 k^2-162\right)-k^4+2 k^2-1=0
\end{equation} 
and has solution $X=dn\left(\frac{1}{3}F[\frac{\pi}{6},k^2],k^2\right)$.  
   
\section{Values of $dn$}

If $r=1$
$$
dn\left(\frac{1}{3}K,m\right)=\sqrt{\frac{1}{2} \left(1+\sqrt{-3+2 \sqrt{3}}\right)}
$$
If $r=2$  
$$
dn\left(\frac{1}{3}K,m\right)=\sqrt{-5+4 \sqrt{2}+\sqrt{51-36 \sqrt{2}}} 
$$ 
If $r=3$
$$
dn\left(\frac{1}{3}K,m\right)=\sqrt{\frac{1}{2}+\frac{\sqrt{2-2^{1/3}}}{2^{5/6}}}
$$
If $r=2/3$
$$
dn\left(\frac{1}{3}K,m\right)=\sqrt{-1 + \sqrt{3}}
$$
If $r=4/3$
$$
dn\left(\frac{1}{3}K,m\right)=\sqrt{2 \left(-2+\sqrt{6}\right)}
$$
If $r=5/3$
$$
dn\left(\frac{1}{3}K,m\right)=\frac{\sqrt{2+\sqrt{3}}}{2}
$$
If $r=7/3$
$$
dn\left(\frac{1}{3}K,m\right)=\sqrt{\frac{1}{2} \left(1+\sqrt{-15+6 \sqrt{7}}\right)}
$$
If $r=8/3$
$$
dn\left(\frac{1}{3}K,m\right)=\sqrt{2 \left(-5-3 \sqrt{3}+\sqrt{57+33 \sqrt{3}}\right)}
$$
If $r=10/3$
$$
dn\left(\frac{1}{3}K,m\right)=\sqrt{-19+8 \sqrt{6}+\sqrt{735-300 \sqrt{6}}}
$$
If $r=11/3$
$$
dn\left(\frac{1}{3}K,m\right)=\sqrt{\frac{26+15 \sqrt{3}+\sqrt{1299+750 \sqrt{3}}}{52+30 \sqrt{3}}}
$$
If $r=13/3$
$$
dn\left(\frac{1}{3}K,m\right)=\frac{1}{2} \sqrt{2+\sqrt{\frac{3}{2} \left(-1+\sqrt{13}\right)}}
$$
If $r=14/3$
$$
dn\left(\frac{1}{3}K,m\right)=\sqrt{-55+21 \sqrt{7}+4 \sqrt{381-144 \sqrt{7}}}
$$
If $r=16/3$
$$
dn\left(\frac{1}{3}K,m\right)=\sqrt{-88-36 \sqrt{6}+6 \sqrt{436+178 \sqrt{6}}}
$$
If $r=17/3$
$$
dn\left(\frac{1}{3}K,m\right)=\sqrt{\frac{1}{2}+\frac{\sqrt{3 \left(1+\sqrt{17}-2 \left(2 \left(1+\sqrt{17}\right)\right)^{1/3}\right)}}{2^{5/6} \left(1+\sqrt{17}\right)^{1/3}}}
$$
If $r=19/3$
$$
dn\left(\frac{1}{3}K,m\right)=\sqrt{\frac{170+39 \sqrt{19}+\sqrt{57459+13182 \sqrt{19}}}{340+78 \sqrt{19}}}
$$
If $r=20/3$
$$
dn\left(\frac{1}{3}K,m\right)=\sqrt{2 \left(-104+33 \sqrt{10}+\sqrt{30 \left(721-228 \sqrt{10}\right)}\right)}
$$
If $r=25/3$
$$
dn\left(\frac{1}{3}K,m\right)=\sqrt{\frac{1}{2} \left(1+10^{1/6} \sqrt{3 \left(-2+10^{1/3}\right)}\right)}
$$
If $r=26/3$
$$
dn\left(\frac{1}{3}K,m\right)=\sqrt{-649+375 \sqrt{3}+4 \sqrt{39 \left(1351-780 \sqrt{3}\right)}}
$$
If $r=29/3$
$$
dn\left(\frac{1}{3}K,m\right)
=\frac{1}{2} \sqrt{2+\sqrt{-6+3 \left(\frac{1}{2} \left(27-5 \sqrt{29}\right)\right)^{1/3}+3 \left(\frac{1}{2} \left(27+5 \sqrt{29}\right)\right)^{1/3}}}
$$
If $r=31/3$
$$
dn\left(\frac{1}{3}K,m\right)=
\sqrt{\frac{1520+273 \sqrt{31}+\sqrt{4617759+829374 \sqrt{31}}}{3040+546 \sqrt{31}}}
$$
If $r=34/3$
$$
dn\left(\frac{1}{3}K,m\right)
=\sqrt{-2449+343 \sqrt{51}+20 \sqrt{6 \left(4999-700 \sqrt{51}\right)}}
$$
If $r=41/3$
$$
dn\left(\frac{1}{3}K,m\right)
=\sqrt{\frac{1}{2}+\frac{2^{1/6} \sqrt{3 \left(1+5 \sqrt{41}-8 \left(2+10 \sqrt{41}\right)^{1/3}\right)}}{\left(1+5 \sqrt{41}\right)^{1/3}}}
$$
If $r=49/3$
$$
dn\left(\frac{1}{3}K,m\right)
=\sqrt{\frac{1}{2} \left(1+3\cdot14^{1/6} \sqrt{-3\cdot2^{1/3}+2\cdot7^{1/3}}\right)}
$$

\[
\]

\centerline{\bf References}\vskip .2in

\noindent

[1]: M. Abramowitz and I.A. Stegun. 'Handbook of Mathematical Functions'. Dover Publications, New York. (1972).

[2]: J.V. Armitage and W.F. Eberlein. 'Elliptic Functions'. Cambridge University Press. (2006)
      
[3]: I.S. Gradshteyn and I.M. Ryzhik. 'Table of Integrals, Series and Products'. Academic Press (1980).

[4]: E.T. Whittaker and G.N. Watson. 'A course on Modern Analysis'. Cambridge U.P. (1927).

\end{document}